\documentclass[12pt,reqno]{amsart}

\usepackage{amsbsy,amsfonts,amsmath,amssymb,amscd,amsthm,mathrsfs,verbatim,calc}
\usepackage{graphicx,epsfig,color}
\usepackage{enumitem}
\usepackage[a4paper,text={6.1in,8in},centering,includefoot,foot=1in]{geometry}
\usepackage[colorlinks=true, linkcolor=blue, citecolor=blue, urlcolor=blue, backref]{hyperref}

\DeclareMathOperator{\initial}{in}
\small\normalsize
\theoremstyle{plain}
\newtheorem{mainthm}{Theorem}

\newtheorem{theorem}{Theorem}[section]
\newtheorem{lemma}[theorem]{Lemma}

\newtheorem{corollary}[theorem]{Corollary}
\newtheorem{definition}[theorem]{Definition}

\theoremstyle{definition}

\newtheorem{remark}[theorem]{Remark}

\numberwithin{equation}{section}
\let\oldmarginpar\marginpar
\renewcommand\marginpar[1]{\-\oldmarginpar[\raggedleft\footnotesize \textcolor{red}{#1}]{\raggedright\footnotesize\textcolor{red}{#1}}}

\begin{document}
\title[Gr\"{o}bner bases for determinantal facet ideals of simplicial complexes]{Gr\"{o}bner bases for determinantal facet ideals of simplicial complexes}

\author[F. Khosh-Ahang]{Fahimeh Khosh-Ahang Ghasr}
\address{{Department of Mathematics, Ilam University,
P.O.Box 69315-516, Ilam, Iran.}}
\email{f.khoshahang@ilam.ac.ir and fahime$_{-}$khosh@yahoo.com}

\begin{abstract}
We provide necessary and sufficient conditions for simplicial complexes whose determinantal facet ideals admit reduced Gr\"{o}bner bases under diagonal term orders. Building on and extending foundational results for binomial edge ideals and determinantal ideals, we introduce two new classes of simplicial complexes---strong closed and poor closed---that generalize the notion of closedness in higher dimensions. Our main theorem offers a unified framework that recovers and refines several known results, including those for unit interval graphs and determinantal ideals of complete graphs. In particular, we correct and generalize prior characterizations of Gr\"{o}bner bases for determinantal facet ideals, establishing radicality for strong closed complexes and providing a new proof for the Gr\"{o}bner basis of maximal minors.

\end{abstract}


%	  05E40   Combinatorial aspects of commutative algebra
%     05E45   Combinatorial aspects of simplicial complexes
%     05C75   Structural characterization of families of graphs
%     05C25   graphs and abstract algebra (groups, rings, ...)
%     05C50   graphs and linear algebra
%     05C40   connectivity

%	  06A11   Algebraic aspects of posets

%     11C20   Matrices , determinants in number theory
%     11C08   polynomials and matrices

%     13P10   Grobner basis other bases for ideals and modules
%     13C14   Cohen-macaulay modules
%	  13H10   Special types (Cohen-Macaulay, Gorenstein, Buchsbaum, etc.)
%	  13D02   Syzygies, resolutions, complexes
%	  13A02   Graded rings
%  	  13F20   Polynomial rings and ideals; rings of integer-valued polynomials
%	  13A18   Valuations and their generalizations
%     13C40   linkage, complete intersection and determinantal ideal

%	  14M25   Toric varieties, Newton polyhedra [See also 52B20]

%	  16S36   Ordinary and skew polynomial rings and semigroup rings


\maketitle
\setcounter{tocdepth}{1}
%\tableofcontents

%<<<<<<<<<<<<<<<<<<<<<<<<<<<<<<<<<<<<<<<<<<<<<<
\section{Introduction and preliminaries}
Interval graphs, defined as intersection graphs of intervals on a line, are a cornerstone of applied mathematics, with applications in scheduling, bioinformatics, and network design. Their algebraic counterparts, such as binomial edge ideals and determinantal facet ideals, bridge graph theory and commutative algebra, enabling the study of combinatorial structures through Gr\"{o}bner bases and homological invariants. In this manuscript we introduce two new concepts,  \emph{strong closed} and \emph{poor closed simplicial complexes}, and provide necessary and sufficient conditions for complexes whose determinantal facet ideals admit quadratic Gr\"{o}bner bases.

Throughout, let $n,d \in \mathbb{N}$,  $H = (V(H), E(H))$ be a simple graph on $n$ vertices and $\mathbb{K}$ a field. Let $\Delta$ be a pure $d$-dimensional simplicial complex (briefly $d$-complex) on $[n]=\{1, \dots, n\}$ whose set of facets is denoted by $\mathcal{F}(\Delta)$. We associate to $H$ the simplicial complex $\Delta_d(H)$ with
$$\mathcal{F}(\Delta_d(H))=\{U \subseteq V(H) \ \mid \text{ the induced subgraph of  } H \text{ on } U \text{ is connected},  |U| = d+1 \};$$ 
for $d=1$, its facet ideal is the edge ideal of $H$.  

Interval graphs and their variants---proper and unit interval graphs---have been linked to binomial edge ideals \cite{B+S+V, Herzog2, Herzog1}. In \cite{Herzog1}, \emph{closed graphs}, equivalent to unit interval graphs, are characterized as those whose binomial edge ideals have quadratic Gr\"{o}bner bases under diagonal term orders. The reference \cite{B+S+V} provides a hierarchy of such generalizations, motivating extensions to simplicial complexes.

The determinantal facet ideal $J_\Delta$ of a $d$-complex $\Delta$ generalizes the binomial edge ideal \cite{Ene}. Efforts to extend \cite[Theorem 1.1]{Herzog1} to simplicial complexes \cite{Ene} faced limitations \cite[Remark 85]{B+S+V}, and attempts to classify complexes with quadratic Gr\"{o}bner bases for $J_\Delta$ \cite{B+S+V} were incomplete. We address these gaps by introducing generalizations of closedness: strong closed and poor closed simplicial complexes (Definition~\ref{closedness}), which captures a broader class of complexes with desirable algebraic properties. These notions arise naturally from a detailed analysis of Buchberger’s criterion and allow us to establish necessary and sufficient conditions for the generators of a determinantal facet ideal to form a reduced Gr\"{o}bner basis.

Despite the computational challenges of large determinants, we establish key structural properties (Lemmas~\ref{determinant} and \ref{determinant2}). These enable our main result, which generalizes \cite[Theorem 1.1]{Herzog1}, refines \cite[Theorem 1.1]{Ene}, and recovers \cite[Theorem 2.7]{Almousa2} and part of \cite[Theorem 82]{B+S+V}. Indeed Theorem \ref{A} provides some necessary and sufficient conditions for complexes with Gr\"{o}bner bases for $J_\Delta$, corrects the ``only if'' part of \cite[Theorem 1.1]{Ene} for $d > 1$, and yields an alternative proof for the Gr\"{o}bner basis of $I_{d+1}(X)$. So we can conclude $J_\Delta$ is a radical ideal for  strong closed $d$-complexes.

This work advances the algebraic combinatorics of higher-dimensional structures, offering new tools for Gr\"{o}bner theory.

\section{Main result}
Let $X=(x_{p,q})_{m\times n}$ be a matrix of indeterminates over $\mathbb{K}$ and $t\leq \min \{m, n\}$ be a positive integer. The ideal $I_t(X)$ generated by all $t$-minors of $X$ in the polynomial ring $S = \mathbb{K}[x_{p,q} \mid 1 \leq p \leq m, 1 \leq q \leq n]$ is called a \emph{determinantal ideal}. These ideals have been studied from various perspectives over the years, connecting commutative algebra with fields such as algebraic geometry, combinatorics, invariant theory, and representation theory (cf. \cite{BC, BV}).

Hereafter $m=d+1$. For each $f,g\in S$ and $1\leq j,j' \leq n$, we set
$$fj-gj'=\begin{pmatrix}
fx_{1,j}-gx_{1,j'} \\
\vdots \\
fx_{d+1,j}-gx_{d+1,j'}
\end{pmatrix}.$$
If $1\leq i_1< \dots < i_t\leq d+1$ and  $C=\{\ell_1, \dots, \ell_t\}\subseteq \{fj-gj' \mid f,g\in S, 1\leq j,j' \leq n\}$, then we use the notation $(C)$ or $(\ell_1, \dots, \ell_t)_{i_1, \dots, i_t}$ for the matrix with respect to rows $i_1, \dots, i_t$ and columns $\ell_1, \dots, \ell_t$ when there is no ambiguity about the order of columns. If moreover $\{\ell_1 , \dots , \ell_t\}\subseteq [n]$, we call $C$ as the \emph{Range} of the submatrix $(C)$ of $X$ and denote it by $\mathrm{Range}(C)$. If $t=d+1$, we set $(\ell_1, \dots, \ell_{d+1})$ instead of $(\ell_1, \dots, \ell_{d+1})_{1, \dots, d+1}$. Moreover, we use the conventional notations $\widehat{\ell}$ for removing $\ell$ from an ordered list and $f\pm g$ for $f+g-g$. We write any set $F=\{i_1, \dots, i_s\}$ with $1\leq i_1< \dots < i_s\leq n$ as $F: i_1 \dots i_s$ where $1\leq s \leq n$. For given two distinct sets $F: i_1 \dots i_{d+1}$ and $G: j_1 \dots j_{d+1}$ we write 
$$\mathcal{A}_{F,G}=\{a \ | \ 1\leq a \leq d+1, i_a\neq j_a\}=a_1 \dots a_r.$$
 Also for any integers $k,b,b'$  with $1\leq k \leq r$, $1\leq b \leq  a_k-1$ and $a_k+1 \leq b' \leq d+1$ we set
 $$S_{k,b}(F,G)=(j_1, \dots, j_{b-1}, \widehat{j_b}, j_{b+1}, \dots  , j_{a_k},i_{a_k}, \dots , i_{d+1}),$$
  and
   $$S_{k,b'}(F,G)=(j_1, \dots , j_{a_k},i_{a_k}, \dots , i_{b'-1}, \widehat{i_{b'}}, i_{b'+1}, \dots  i_{d+1}).$$ 
Note that $S_{k,b}(F,G)$ and $S_{k,b'}(F,G)$ may have some equal columns. 

Let $\Delta$ be a pure simplicial complex on $[n]$ of positive dimension $d$. There are various attempts to generalize the concept of binomial edge ideals (cf. \cite{Almousa, Ene, M+R}). One of them is the \emph{determinantal facet ideal} $J_\Delta$ of $\Delta$ introduced in \cite{Ene}:
$$J_\Delta=\langle \det(i_1, \dots, i_{d+1}) \mid i_1 \dots i_{d+1} \in \mathcal{F}(\Delta)\rangle.$$
Clearly, $J_{\Delta_1(H)}$ is the binomial edge ideal $J_H$ of $H$ introduced in \cite{Herzog1}. Also, if $H$ is a complete graph, then $J_{\Delta_d(H)}$ is the determinantal ideal $I_{d+1}(X)$ if $d+1\leq n$.

The concept of closed graphs can be generalized from different perspectives (cf. \cite{Almousa, B+S+V, Ene}). For the aim of this note, we use the following four generalizations:
\begin{definition}\label{closedness}
\begin{enumerate}
    \item (\cite{Ene}) $\Delta$ is called a \emph{closed simplicial complex} if there is a labelling $[n]$ on $V(\Delta)$ such that for every two distinct facets $F: i_1\dots i_{d+1}$ and $G: j_1 \dots j_{d+1}$ when $i_t=j_t$ for some $1\leq t \leq d+1$, $\Delta$ contains the full $d$-skeleton of the simplex on the vertex set $F\cup G$.
    \item (\cite[Definition 30]{B+S+V}) $\Delta$ is called a \emph{unit interval simplicial complex} if there is a labelling $[n]$ on $V(\Delta)$ such that for each facet $F: i_1\dots i_{d+1}$, $\Delta$ contains the full $d$-skeleton of the simplex on the vertex set $\{i_1, i_1+1,i_1+2, \dots , i_{d+1}\}$.
    \item $\Delta$ is called a \emph{poor closed simplicial complex} if there is a labelling $[n]$ on $V(\Delta)$ such that for every two distinct facets $F: i_1 \dots i_{d+1}$ and $G: j_1 \dots j_{d+1}$ when $i_t=j_t$ for some $1\leq t \leq d+1$, there exists a facet contained in $F\cup G$ other than $F$ and $G$.
    \item $\Delta$ is called a \emph{strong closed simplicial complex} if there is a labelling $[n]$ on $V(\Delta)$ such that for every two distinct facets $F: i_1 \dots i_{d+1}$ and $G: j_1 \dots j_{d+1}$   when $i_1\leq j_1$ and $\mathcal{A}_{F,G}=a_1\dots a_r$ for some $r$ with $1\leq r \leq d$, $\Delta$ contains $\mathrm{Range}S_{k,b}(F,G)$ and $\mathrm{Range}S_{k,b'}(F,G)$  for any integers $k,b,b'$ with $1\leq k \leq r, 1\leq b \leq  a_k-1$ and 
        $a_k+1 \leq b' \leq d+1$.
\end{enumerate}
% $H$ is called $d$-unit interval (resp. $d$-closed, $d$-poor closed, $d$-strong closed) if $\Delta_d(H)$ is so.
\end{definition}
We start by finding the relations between these concepts with the next lemma and remarks.
\begin{lemma}\label{lem1}
Every closed or unit interval simplicial complex is strong closed and every strong closed simplicial complex is poor closed with the same labelling on $V(\Delta)$.
\end{lemma}
\begin{proof}
Clearly each closed simplicial complex is strongly closed. Suppose $\Delta$ is unit interval, $F:i_1 \dots i_{d+1}$ and $G:j_1\dots j_{d+1}$ are two facets with $i_1\leq j_1$ and $\mathcal{A}_{F,G}=a_1\dots a_r$ for some $1\leq r \leq d$. Let $k,b,b'$ be integers with $1\leq k \leq r$,  $1\leq b \leq a_k-1$,  $a_k+1\leq b'\leq d+1$.
If $j_{d+1}\leq i_{d+1}$, then $i_1\leq j_1, \dots , j_{a_k}, i_{a_k}, \dots , i_{d+1} \leq i_{d+1}$. So since $F\in \Delta$ and $\Delta$ is unit interval, $\Delta$ contains $\mathrm{Range}S_{k,b}(F,G)$ and $\mathrm{Range}S_{k,b'}(F,G)$. Now suppose $i_{d+1}<j_{d+1}$. If
 $$\{\min \{j_1, \dots , j_{a_k}, i_{a_k}, \dots , i_{d+1}\} , \max \{j_1, \dots , j_{a_k}, i_{a_k}, \dots , i_{d+1}\}\}$$
 is contained in $F$ (resp. $G$), then $\Delta$ contains $\mathrm{Range}S_{k,b}(F,G)$ and $\mathrm{Range}S_{k,b'}(F,G)$ by unit intervalness of $\Delta$. Otherwise since $r\leq d$, 
 $$\min  \{j_1, \dots , j_{a_k}, i_{a_k}, \dots , i_{d+1}\} =j_{a_1}=j_1, \ \max  \{j_1, \dots , j_{a_k}, i_{a_k}, \dots , i_{d+1}\} =i_{a_r}=i_{d+1} .$$
 Thus since $G\in \Delta $,    $ j_1\leq j_1, \dots , j_{a_k}, i_{a_k}, \dots , i_{d+1}\leq j_{d+1}$ and $\Delta$ is unit interval, $\Delta$ contains $\mathrm{Range}S_{k,b}(F,G)$ and $\mathrm{Range}S_{k,b'}(F,G)$. These prove the first statement.

Suppose $\Delta$ is strong closed, $F:i_1 \dots i_{d+1}$ and $G:j_1\dots j_{d+1}$ are two facets with $i_1\leq j_1$ and $\mathcal{A}_{F,G}=a_1\dots a_r$ for some $1\leq r \leq d$. If $j_{a_1}\neq i_b$ for all $a_1+1\leq b \leq d+1$, then $\mathrm{Range}S_{1,b}(F,G)$ is a facet except $F$ and $G$ for each $b$ with $a_1+1 \leq b \leq d+1$. Otherwise $j_{a_1}=i_b$ for some $b$ with $a_1+1\leq b \leq d+1$. In this case $\mathrm{Range}S_{1,b}(F,G)$ is a facet other than $F$ and $G$.
\end{proof}

\begin{remark}\label{remark}
\begin{enumerate}
\item If for some labeling $[n]$ on $V(\Delta)$ any two facets  $F: i_1 \dots i_{d+1}$ and $G: j_1 \dots j_{d+1}$ in $\Delta$ satisfies $i_t\neq j_t$ for all $1\leq t \leq d+1$, then clearly $\Delta$ is closed, poor closed and strong closed.
\item The simplicial complex $\Delta$ with
$$\mathcal{F}(\Delta)=\{123, 124, 134, 234, 235, 245, 345, 568, 789, 8\ 10\ 11\},$$
    introduced in \cite[Remark 84]{B+S+V},  is neither unit interval nor closed while it is strong closed and so poor closed. This example, in conjunction with Theorem \ref{A}, justifies \cite[Remark 84]{B+S+V}. 
\item The class of poor closed simplicial complexes properly contains the class of strong closed simplicial complexes. For instance, the simplicial complex $\Delta$ with 
$$\mathcal{F}(\Delta)=\{123, 124, 234\}$$
is a poor closed $2$-complex which is not strong closed.
\item In view of Parts 2 and 3 and Lemma \ref{lem1}, 
 $$\{\text{ unit interval } \}\cup \{ \text{ closed } \} \subset \{ \text{ strong closed } \} \subset \{ \text{ poor closed }\}.$$
 Notice Lemmas 42 and 43 in \cite{B+S+V} indicate there is no inclusion relation between the classes of closed simplicial complexes and unit interval simplicial complexes.
\item For a connected graph $H$ and $\Delta=\Delta_1(H)$, the concepts of unit interval, closed, strong closed and poor closed simplicial complexes are equivalent to the concepts of unit interval graphs and closed graphs, but this is not true for $d$-complexes when $d>1$ as mentioned in (4).
\item Every unit interval simplicial complex $\Delta$ is a flag complex and thus a clique complex. Hence, there exists a graph $H$ whose maximal cliques are exactly the facets of $\Delta$. So, $H$ is a proper interval graph. Therefore, each unit interval simplicial complex is a clique complex of a proper interval graph.
\end{enumerate}
\end{remark}

The following two lemmas not only play  essential roles in the proof of the main result, but also provide  special properties of determinant.

\begin{lemma}\label{determinant}
For distinct integers $\ell _1, \dots , \ell _d , j_1$ in $[n]$ we have
\begin{align*}\det(x_{1,j_1}\ell _1&-x_{1,\ell _1}j_1, \ell _2 , \dots , \ell _d)_{2, \dots , d+1}=\\
&\det(j_1, \ell _1, \dots , \ell _d)+\sum_{2\leq r\leq d} (-1)^{r+1}x_{1,\ell _r}\det(j_1, \ell _1,\dots , \widehat{\ell _r},\dots ,  \ell _d)_{2,\dots , d+1}.
\end{align*}
\end{lemma}
\begin{proof}
We proceed by induction on $d$. The result is clear for $d=1$. Assume inductively $d>1$ and the result has been proved for smaller values of $d$. Applying Laplace expansion along the last column, the inductive hypothesis for the gained minors, adding and subtracting the same needed value, and finally using Laplace expansion again imply the result as follows:
\begin{align*}
\det(&x_{1,j_1}\ell _1-x_{1,\ell _1}j_1, \ell _2 , \dots , \ell _d)_{2, \dots , d+1}\\
&=\sum_{2\leq k \leq d+1} (-1)^{(k -1)+d}x_{k ,\ell _d}\det(x_{1,j_1}\ell _1-x_{1,\ell _1}j_1, \ell _2 , \dots , \ell _{d-1})_{2, \dots, \widehat{k}, \dots , d+1}\\
&=\sum_{2\leq k \leq d+1} (-1)^{(k-1)+d}x_{k,\ell _d} \Bigg( \det(j_1, \ell _1, \dots , \ell _{d-1})_{1, \dots , \widehat{k}, \dots , d+1}\\
&\hspace{4mm}+\sum_{2\leq r\leq d-1} (-1)^{r+1}x_{1,\ell _r}\det(j_1, \ell _1, \dots, \widehat{\ell _r} , \dots , \ell _{d-1})_{2, \dots, \widehat{k} , \dots , d+1}
\Bigg)\\
&= \Bigg ( \sum_{2\leq k \leq d+1} (-1)^{k+(d+1)}x_{k,\ell _d}\det(j_1, \ell _1, \dots , \ell _{d-1})_{1 , \dots , \widehat{k} , \dots , d+1}\\
&\hspace{4mm}\pm (-1)^{1+(d+1)}x_{1,\ell _d}\det(j_1, \ell _1, \dots , \ell _{d-1})_{2, \dots , d+1}\Bigg ) \\
&\hspace{4mm}+ \sum_{2\leq r\leq d-1} (-1)^{r+1}x_{1,\ell _r}\Bigg (\sum_{2\leq k \leq d+1} (-1)^{(k-1)+d}x_{k,\ell _d}
\det(j_1, \ell _1, \dots , \widehat{\ell _r}, \dots , \ell _{d-1})_{2, \dots , \widehat{k} , \dots , d+1}\Bigg )\\
&= \Bigg(\det(j_1, \ell _1, \dots , \ell _d)- (-1)^{1+(d+1)}x_{1,\ell _d}\det(j_1, \ell _1, \dots , \ell _{d-1})_{2, \dots , d+1}\Bigg)\\
&\hspace{4mm} + \sum_{2\leq r\leq d-1} (-1)^{r+1}x_{1,\ell _r}\det(j_1, \ell _1, \dots , \widehat{\ell _r} , \dots , \ell _d)_{2, \dots , d+1}\\
&=\det(j_1, \ell _1, \dots , \ell _d)+\sum_{2\leq r\leq d} (-1)^{r+1}x_{1,\ell _r}\det(j_1, \ell _1,\dots , \widehat{\ell _r},\dots ,  \ell _d)_{2,\dots , d+1}.
\end{align*}
\end{proof}

\begin{lemma}\label{determinant2}
Suppose $F:i_1\dots i_{d+1}$ and $G:j_1, \dots , j_{d+1}$ with $\mathcal{A}_{F,G}=a$. Then
$$\det(i_1,\dots,i_{a-1}, x_{a,j_a}i_a - x_{a,i_a}j_a, i_{a+1},\dots,i_{d+1})
= \sum_{\substack{1\leq k \leq d+1 \\ k\neq a}} (-1)^{1+k} x_{a,i_k}
\det S_{1,k}(F,G).$$
\end{lemma}
\begin{proof}
Since the determinant’s value is unchanged (up to sign) with row and column swaps, assume $a=1$. Suppose $M_{1,i_k}$ is the minor of the element $x_{1,i_k}$  for $k$ with $2\leq k \leq d+1$ in
 $$(x_{1,j_1}i_1 - x_{1,i_1}j_1, i_2, \dots, i_{d+1}).$$
 Applying Lemma \ref{determinant} yields
\begin{align*}
\det(x_{1,j_1}i_1 &- x_{1,i_1}j_1, i_2, \dots, i_{d+1}) \\
&=\sum_{2\leq k \leq d+1} (-1)^{1+k} x_{1,i_k}M_{1,i_k}\\
 &=\sum_{2\leq k \leq d+1} (-1)^{1+k} x_{1,i_k} \det(x_{1,j_1}i_1-x_{1,i_1}j_1, i_2,\dots, \widehat{i_k}, \dots , i_{d+1})_{2, \dots , d+1}\\
 &=\sum_{2\leq k \leq d+1} (-1)^{1+k} x_{1,i_k} \Bigg(
\det(j_1, i_1, \dots , \widehat{i_k}, \dots , i_{d+1}) \\
&\hspace{4mm} +\sum_{2\leq r\leq k-1} (-1)^{r+1}x_{1,i_r}\det(j_1, i_1,\dots , \widehat{i_r}, \dots, \widehat{i_k}, \dots ,  i_{d+1})_{2,\dots , d+1}\\
 &\hspace{4mm} +\sum_{k+1\leq r\leq d+1} (-1)^rx_{1,i_r}\det(j_1, i_1,\dots , \widehat{i_k},\dots, \widehat{i_r}, \dots ,  i_{d+1})_{2,\dots , d+1}\Bigg)\\
 &=\sum_{2\leq k \leq d+1} (-1)^{1+k} x_{1,i_k} \det(j_1, i_1, \dots , \widehat{i_k}, \dots , i_{d+1}) \\
&\hspace{4mm} +\sum_{2\leq k \leq d+1}\Bigg(\sum_{2\leq r\leq k-1}  (-1)^{k+r} x_{1,i_k}x_{1,i_r}\det(j_1, i_1,\dots , \widehat{i_r}, \dots, \widehat{i_k}, \dots ,  i_{d+1})_{2,\dots , d+1}\\
 &\hspace{4mm} +\sum_{k+1\leq r\leq d+1} (-1)^{k+r+1} x_{1,i_k}x_{1,i_r}\det(j_1, i_1,\dots , \widehat{i_k},\dots, \widehat{i_r}, \dots ,  i_{d+1})_{2,\dots , d+1}\Bigg)\\
 &=\sum_{2\leq k \leq d+1} (-1)^{1+k} x_{1,i_k} \det(j_1, i_1, \dots , \widehat{i_k}, \dots , i_{d+1}).
\end{align*}
\end{proof}

Hereafter, we consider the lexicographic order induced by $x_{1,1}>\dots >x_{1,n}>x_{2,1}>\dots >x_{2,n}> \dots > x_{d+1,1}> \dots > x_{d+1,n}$. As known, the set of all $(d+1)$-minors of $X$ is a Gr\"{o}bner basis for $I_{d+1}(X)$ (\cite{B+Z, S+Z}). Also, in \cite[Theorem 1.1]{Ene}, it is shown that $\{\det(j_1, \dots, j_{d+1}) \mid j_1 \dots j_{d+1}\in \mathcal{F}(\Delta)\}$ is a Gr\"{o}bner basis for $J_\Delta$ if and only if $\Delta$ is closed. However, according to \cite[Remark 85]{B+S+V}, the \emph{only if} part of this result is not correct. Although the \emph{if} part of \cite[Theorem 1.1]{Ene} is clear using Buchberger's criterion, we prove it again and show explicitly how the desired $S$-pairs reduce to zero. This helps us generalize \cite[Theorem 1.1]{Herzog1}, improve \cite[Theorem 1.1]{Ene}, and regain \cite[Theorem 2.7]{Almousa2} and part of Theorem 82 in \cite{B+S+V}.  This also provides an alternative, more understandable proof that the set of all $(d+1)$-minors of $X$ is a Gr\"{o}bner basis for $I_{d+1}(X)$ and illustrates that the \emph{only if} part of \cite[Theorem 1.1]{Ene} cannot hold for $d>1$, since strong closed simplicial complexes are not necessarily closed. 

\begin{mainthm}\label{A}
Suppose $\Delta$ is a $d$-dimensional simplicial complex on $[n]$ and
$$\mathcal{B}=\{\det(i_1, \dots , i_{d+1}) \ | \ i_1 \dots  i_{d+1}\in \mathcal{F}(\Delta )\}.$$
\begin{enumerate}
  \item The set of all $(d+1)$-minors of $X=(x_{p,q})_{(d+1)\times n}$ forms a reduced Gr\"{o}bner basis for $I_{d+1}(X)$.
  \item  If $\Delta$ is a strong closed simplicial complex, then  $\mathcal{B}$ forms a reduced Gr\"{o}bner basis for $J_\Delta $.
  \item If $\mathcal{B}$ forms a Gr\"{o}bner basis for $J_\Delta$, then $\Delta$ is a poor closed simplicial complex.
  \item $H$ is a closed graph if and only if $\{\det(i,j) \ | \ \{i,j\}\in E(H)\}$ forms a reduced Gr\"{o}bner basis for $J_H$.
\end{enumerate}
\end{mainthm}
\begin{proof}
(1) We use Buchberger's criterion. Suppose $F: i_1 \dots  i_{d+1}, G: j_1\dots  j_{d+1}$ and $f=\det(i_1, \dots, i_{d+1}), g=\det(j_1, \dots, j_{d+1})$. We show $S(f,g)$ reduces to zero in the following three cases:

    \textbf{Case I.} Let $\mathcal{A}_{F,G}=1 \dots d+1$. Then the initial monomials of $f$ and $g$ are relatively prime, so $S(f,g)$ reduces to zero.
% by \cite[Lemma 2.3.1]{H+H}.

    \textbf{Case II.} Let $\mathcal{A}_{F,G}=a$ and  $i_a<j_a$. Then by elementary properties of determinants together with Lemma \ref{determinant2} we obtain
\begin{align*}
S(f,g)&=x_{a,j_a}\det(i_1, \dots , i_{d+1})-x_{a,i_a}\det(j_1, \dots , j_{d+1})\\
&=\det(i_1, \dots , i_{a-1}, x_{a,j_a}i_a, i_{a+1}, \dots , i_{d+1})-\det(j_1, \dots , j_{a-1}, x_{a,i_a}j_a,j_{a+1}, \dots , j_{d+1}) \\
&=\det(i_1,\dots , i_{a-1}, x_{a,j_a}i_a-x_{a,i_a}j_a, i_{a+1}, \dots,i_{d+1}) \\
&=\sum_{\substack{1\leq k \leq d+1\\ k\neq a}} (-1)^{1+k} x_{a,i_k} \det S_{1,k}(F,G).
\end{align*}
Now it remains to compare the initial terms of $S(f,g)$ and $x_{a,i_k} \det S_{1,k}(F,G)$. We know the initial term of $S(f,g)$ arises from $x_{a,i_a}\det(j_1,\dots,j_{d+1})$  or $x_{a,j_a}\det(i_1, \dots , i_{d+1})$  which is not cancelled in the subtraction. As the initial monomial is determined by the product along the main diagonal, whenever the choice of the main diagonal would eliminate a term, in a $2\times 2$ minor the corresponding secondary diagonal must be taken instead. So since $i_a<j_a$, the initial term of $S(f,g)$ arises from $x_{a,i_a}\det(j_1,\dots,j_{d+1})$ except when $a=d+1$. Therefore
$$ \initial (S(f,g)) = 
\begin{cases}
    x_{a,i_a}x_{a, i_{a+1}}x_{a+1, j_a}\prod_{\substack{1\leq \ell \leq d+1 \\ \ell\neq a,a+1}} x_{\ell,i_\ell}, & \text{if } 1\leq a\leq d, \\
    x_{d+1,j_{d+1}}x_{d, i_{d+1}}x_{d+1, i_d}\prod_{\substack{1\leq \ell \leq d-1}} x_{\ell,i_\ell}, & \text{if } a=d+1,
\end{cases}$$
and
\begin{align*}
&\mathrm{in}(x_{a,i_k} \det S_{1,k}(F,G))\\
&=\begin{cases}
\mathrm{in}(x_{a,i_k} \det(i_1, \dots , \widehat{i_k}, \dots i_{a-1}, i_a, j_a, i_{a+1} \dots , i_{d+1})), & \text{if } 1\leq k \leq a-1, \\
\mathrm{in}(x_{a,i_k} \det(i_1, \dots , i_a, j_a, i_{a+1} \dots \widehat{i_k}, \dots , i_{d+1})), & \text{if } a+1\leq k \leq d+1,
\end{cases}\\
&=\begin{cases}
x_{a,i_k} (\prod_{1\leq \ell \leq k-1}x_{\ell,i_\ell})  (\prod_{k\leq \ell \leq a-1}x_{\ell,i_{\ell+1}}) x_{a,j_a} (\prod_{a+1\leq \ell \leq d+1}x_{\ell,i_\ell}) , & \text{if } 1\leq k \leq a-1, \\
x_{a,i_k} (\prod_{1\leq \ell \leq a}x_{\ell,i_\ell}) x_{a+1,j_a} (\prod_{a+2\leq \ell \leq k} x_{\ell, i_{\ell-1}}) (\prod_{k+1\leq \ell \leq d+1}x_{\ell, i_\ell}), & \text{if } a+1\leq k \leq d+1.
\end{cases}\\
\end{align*}
Therefore if either $k=a+1$ or $a=d+1, k=d$, then two initial terms equal.  Else when $1\leq k \leq a-1$  the first variables which differ between $\initial (S(f,g))$ and $\mathrm{in}(x_{a,i_k} \det S_{1,k}(F,G))$ are from the $k$-th row and are respectively $x_{k,i_k}$ and $x_{k, i_{k+1}}$. Also when $a+2 \leq k \leq d+1$ the first variables which differ between $\initial (S(f,g))$ and $\mathrm{in}(x_{a,i_k} \det S_{1,k}(F,G))$ are from the $a$-th row and are respectively $x_{a,i_a}x_{a,i_{a+1}}$ and $x_{a,i_k}x_{a, i_a}$. These observations indicate $\mathrm{in}(S(f,g)) \geq \mathrm{in}(x_{1,i_k} \det S_{1,k}(F,G))$. Thus, $S(f,g)$ reduces to zero.
 
 \textbf{Case III.} 
Let $\mathcal{A}_{F,G}=a_1\dots a_r$ for some $r$ with $2\leq r \leq d$ and $i_{a_1}<j_{a_1}$. Elementary properties of determinants, adding and subtracting the same needed values, and at last  applying Lemma \ref{determinant2} yield
\begin{align*}
&S(f,g)\\
=&(\prod_{a\in \mathcal{A}_{F,G}}x_{a,j_a})\det(i_1, \dots , i_{d+1})-(\prod_{a\in \mathcal{A}_{F,G}}x_{a,i_a})\det(j_1, \dots , j_{d+1})\\
=&\det(\{i_a \ | \ a\notin \mathcal{A}_{F,G}\}\cup \{x_{a,j_a}i_a \ | \ a\in \mathcal{A}_{F,G}\})-\det(\{j_a \ | \ a\notin \mathcal{A}_{F,G}\}\cup \{x_{a,i_a}j_a \ | \ a\in \mathcal{A}_{F,G}\})\\
=&\det(\{i_a \ | \ i_a=j_a\}\cup \{x_{a,j_a}i_a \ | \ a\in a_1\dots a_r\})\\
 &\pm\sum_{2\leq k \leq r} \det(\{i_a \ | \ i_a=j_a\}\cup \{x_{a,j_a}i_a \ | \ a\in a_k \dots a_r\} \cup \{x_{a,i_a}j_a \ | \ a\in a_1 \dots a_{k-1}\})\\
 &-\det(\{j_a \ | \ i_a=j_a\}\cup \{x_{a,i_a}j_a \ | \ a\in a_1\dots a_r\})\\
 =&\sum_{1\leq k \leq r}\Bigg(\det(\{i_a \ | \ i_a=j_a \}\cup \{x_{a,j_a}i_a \ | \ a\in a_k \dots a_r\} \cup \{x_{a,i_a}j_a \ | \ a\in a_1 \dots  a_{k-1}\})\\
 &-\det(\{i_a \ | \ i_a=j_a \}\cup \{x_{a,j_a}i_a \ | \ a\in a_{k+1} \dots a_r\} \cup \{x_{a,i_a}j_a \ | \ a\in a_1 \dots  a_k\})\Bigg)\\
 =&\sum_{1\leq k \leq r}\det(\{i_a \ | \ i_a=j_a \}\cup \{x_{a,j_a}i_a \ | \ a\in a_{k+1} \dots a_r\} \cup \{x_{a,i_a}j_a \ | \ a\in a_1 \dots  a_{k-1}\}\\
 &\cup \{x_{a_k,j_{a_k}}i_{a_k}-x_{a_k,i_{a_k}}j_{a_k}\})\\
 \end{align*}
 \begin{align*}
 =&\sum_{1\leq k \leq r}\Bigg((\prod _{a\in a_1 \dots  a_{k-1}}x_{a,i_a}) (\prod_{a\in a_{k+1} \dots a_r}x_{a,j_a})\det(j_1, \dots , j_{a_k-1}, x_{a_k,j_{a_k}}i_{a_k}-x_{a_k,i_{a_k}}j_{a_k}, i_{a_k+1}, \dots , i_{d+1})\Bigg)\\
=&\sum_{1\leq k \leq r}\Bigg((\prod _{a\in a_1 \dots  a_{k-1}}x_{a,i_a}) (\prod_{a\in a_{k+1} \dots a_r}x_{a,j_a})
\sum_{\substack{1\leq b \leq d+1\\ b\neq a_k}}(-1)^{1+b}x_{a_k,\ell_b}\det S_{k,b}(F,G)\Bigg)\\
=&\sum_{1\leq k \leq r}\ \ \sum_{\substack{1\leq b \leq d+1\\ b\neq a_k}}\Bigg((-1)^{1+b}x_{a_k,\ell_b}(\prod _{a\in a_1 \dots  a_{k-1}}x_{a,i_a}) (\prod_{a\in a_{k+1} \dots a_r}x_{a,j_a})
\det S_{k,b}(F,G) \Bigg),
\end{align*}
where for each $b$ with $1\leq b\leq a_k-1$, $\ell_b=j_b$ and for each $b$ with $a_k+1\leq b \leq d+1$, $\ell_b=i_b$. 

To compare the initial terms, note that according to above calculations we have 
\begin{align*}
&S(f,g) =\\
&\sum_{1\leq k \leq r}\Bigg((\prod _{a\in a_1 \dots  a_{k-1}}x_{a,i_a}) (\prod_{a\in a_{k+1} \dots a_r}x_{a,j_a})\det(j_1, \dots , j_{a_k-1}, x_{a_k,j_{a_k}}i_{a_k}-x_{a_k,i_{a_k}}j_{a_k}, i_{a_k+1}, \dots , i_{d+1})\Bigg).
\end{align*}
On the other hand, analogous argument to Case II, shows that the initial  term \linebreak $(\prod _{a\in a_1 \dots  a_{k-1}}x_{a,i_a}) (\prod_{a\in a_{k+1} \dots a_r}x_{a,j_a})\det(j_1, \dots , j_{a_k-1}, x_{a_k,j_{a_k}}i_{a_k}-x_{a_k,i_{a_k}}j_{a_k}, i_{a_k+1}, \dots , i_{d+1})$ depends on $k$, since the variables in $a_k$-th row is different. Hence
\begin{align*}
\mathrm{in}(S(f,g)) &= \\
&\max  \{\mathrm{in}\Bigg( (\prod _{a\in a_1 \dots  a_{k-1}}x_{a,i_a}) (\prod_{a\in a_{k+1} \dots a_r}x_{a,j_a})\\
&\det(j_1, \dots , j_{a_k-1}, x_{a_k,j_{a_k}}i_{a_k}-x_{a_k,i_{a_k}}j_{a_k}, i_{a_k+1}, \dots , i_{d+1})\Bigg) \ | \ 1\leq k \leq r\}.
\end{align*}
So by using similar argument to Case II we have 
\begin{align*}
&\mathrm{in}(S(f,g)) \\
& \geq  (\prod _{a\in a_1 \dots  a_{k-1}}x_{a,i_a}) (\prod_{a\in a_{k+1} \dots a_r}x_{a,j_a})\mathrm{in}(\det(j_1, \dots , j_{a_k-1}, x_{a_k,j_{a_k}}i_{a_k}-x_{a_k,i_{a_k}}j_{a_k}, i_{a_k+1}, \dots , i_{d+1}))\\
& \geq (\prod _{a\in a_1 \dots  a_{k-1}}x_{a,i_a}) (\prod_{a\in a_{k+1} \dots a_r}x_{a,j_a}) \mathrm{in} ( x_{a_k,\ell_b} \det S_{k,b}(F,G))\\
&=  \mathrm{in} ( x_{a_k,\ell_b}(\prod _{a\in a_1 \dots  a_{k-1}}x_{a,i_a}) (\prod_{a\in a_{k+1} \dots a_r}x_{a,j_a})
\det S_{k,b}(F,G)),
\end{align*}
for each $k,b$ with $1\leq k \leq r$ and $1\leq b \leq d+1, b\neq a_k$. Hence $S(f,g)$ reduces to zero.

To prove the Gr\"{o}bner basis is reduced, suppose $f=\det(i_1, \dots , i_{d+1})$ and $g=\det(j_1, \dots , j_{d+1})$, where $1\leq i_1< \dots < i_{d+1}\leq n$ and  $1\leq j_1< \dots < j_{d+1}\leq n$. Then $\mathrm{in}(f)=\prod_{1\leq \ell \leq d+1}x_{\ell,i_\ell}$ and each monomial in $\mathrm{supp}(g)$ is of the form $\prod_{1\leq \ell \leq d+1}x_{\ell,\sigma (j_\ell)}$ for some permutation $\sigma$ on $\{1, \dots , d+1\}$. The coefficient of $\mathrm{in}(f)$ is $1$ and if  $\prod_{1\leq \ell \leq d+1}x_{\ell,i_\ell} \ \mid \ \prod_{1\leq \ell \leq d+1}x_{\ell,\sigma (j_\ell)}$, then $\sigma$ is the identity, so $f=g$. This completes the proof.

(2) Consider two facets $F: i_1 \dots  i_{d+1}$ and  $G:  j_1 \dots j_{d+1}$ with $i_1\leq j_1$ and apply the proof of Part (1). Note that if $S_{k,b}(F,G)$ has equal columns, then $\det S_{k,b}(F,G)=0$ and otherwise $\mathrm{Range}S_{k,b}(F,G)$ is a facet of $\Delta$. So  $\det S_{k,b}(F,G)\in \mathcal{B}$ as desired. 

(3) Suppose in contrary there exist distinct facets $F: i_1 \dots  i_{d+1}$ and  $G: j_1 \dots  j_{d+1}$ such that $i_t=j_t$ for some $1\leq t \leq d+1$, but $\Delta$ contains no facet in $F\cup G$ except $F$ and $G$. Since $f=\det(i_1, \dots , i_{d+1})$ and $g=\det(j_1, \dots , j_{d+1})$ are in $J_{\Delta}$, $S(f,g)\in J_{\Delta}$. Thus $\mathrm{in}(S(f,g))\in \mathrm{in}(J_{\Delta})$, so $\mathrm{in}(\det(\ell_1, \dots , \ell_{d+1})) \mid \mathrm{in}(S(f,g))$ for some $\ell_1 \dots  \ell_{d+1}\in \mathcal{F}(\Delta )$. On the other hand by similar discussion to Case II of Part 1, if  $\mathcal{A}_{F,G}=a_1\dots a_r$ for some $r$ with $1\leq r \leq d$, we have 
\begin{align*} 
&\mathrm{in} (S(f,g)) = \\
&\begin{cases}
   x_{a_r, j_{a_r+1}} x_{a_r+1, j_{a_r}}(\prod _{a\in \mathcal{A}_{F,G}} x_{a,i_a}) (\prod_{\substack{1\leq \ell \leq d+1 \\ \ell\neq a_r,a_r+1}} x_{\ell,j_\ell}), & \text{if } 1\leq a_r\leq d, i_{a_r}<j_{a_r}, \\
    x_{d,i_{d+1}}x_{d+1, i_d}(\prod _{a\in \mathcal{A}_{F,G}} x_{a,j_a})(\prod_{1\leq \ell \leq d-1} x_{\ell,i_\ell}),  &\text{if } a_r=d+1, i_{a_r}<j_{a_r}, (r=1 \text{ or } a_{r-1}<d), \\
     x_{d,j_{d+1}}x_{d+1, j_d}(\prod _{a\in \mathcal{A}_{F,G}} x_{a,i_a})(\prod_{1\leq \ell \leq d-1} x_{\ell,j_\ell}), & \text{if } a_r=d+1,a_{r-1}=d, i_d<j_d. \\
\end{cases}
\end{align*}
So if $\prod_{1\leq k \leq d+1}x_{k,\ell_k}$ divides $\mathrm{in}(S(f,g))$, then $\{\ell_1, \dots , \ell_{d+1}\}\subseteq F\cup G$ with nonempty intersection with both $F$ and $G$, a contradiction. Note $\prod_{1\leq k \leq d+1}x_{k,\ell_k}$ can't contain the elements of a secondary diagonal of some $2\times 2$ minor  simultaneously, for instance $x_{a_r,j_{a_r+1}}$ and $x_{a_r+1, i_{a_r}}$, since $a_r<a_r+1$ and $j_{a_r+1}>j_{a_r}$.

(4) The result follows from (2) and (3) for $\Delta=\Delta_1(H)$ and  Remark \ref{remark}(5).
\end{proof}

We end this note by the following corollary of Theorem \ref{A} which can be gained by repeating the process used in \cite[Corollary 2.2]{Herzog1}. 

\begin{corollary}
For a strong closed simplicial complex $\Delta$, $J_\Delta$ is a radical ideal.
\end{corollary}

%\textbf{Statements and Declarations.}

%The author declares that no funds, grants, or other support were received during the preparation of this manuscript. Also the author has no relevant financial or non-financial interests to disclose.


\begin{thebibliography}{99}
\bibitem{Almousa2}  A. Almousa K. N. Lin, and W. Liske, {\it Rees algebras of unit interval determinantal facet ideals.} arxiv.2008.10950.

\bibitem{Almousa} A. Almousa and K. VandeBogert, {\it Determinantal facet ideals for smaller minors.} Arch. Math. (Basel) 118 (2022), no. 3, 247--256.

\bibitem{B+S+V} B. Benedetti, L. Seccia and M. Varbaro, {\it Hamiltonian paths, unit-interval complexes, and determinantal facet ideals.} Adv. in Appl. Math. 141 (2022), Paper No. 102407, 55 pp.

\bibitem{B+Z} D. Bernstein and A. Zelevinsky, {\it Combinatorics of maximal minors.} J. Algebraic Combin. 2 (1993), no. 2, 111--121.

\bibitem{BC} W. Bruns and A. Conca, {\it Gröbner bases and determinantal ideals. Commutative algebra, singularities and computer algebra} (Sinaia, 2002), 9--66, NATO Sci. Ser. II Math. Phys. Chem., 115, Kluwer Acad. Publ., Dordrecht, 2003.

\bibitem{BV} W. Bruns and U. Vetter, Determinantal rings. Lecture Notes in Mathematics, 1327. Springer-Verlag, Berlin, 1988. {\rm viii}+236 pp.

\bibitem{Ene} V. Ene, J. Herzog, T. Hibi and F. Mohammadi, {\it Determinantal facet ideals.} Michigan Math. J. 62 (2013), no. 1, 39--57.

\bibitem{Herzog2} V. Ene, J. Herzog and T. Hibi, {\it Cohen-Macaulay binomial edge ideals.} Nagoya Math. J. 204 (2011), 57--68.

\bibitem{Herzog1} J. Herzog, T. Hibi, F. Hreinsd\'{o}ttir, T. Kahle and J. Rauh, {\it Binomial edge ideals and conditional independence statements.}  Adv. in Appl. Math. 45 (2010), no. 3, 317--333.

\bibitem{M+R} F. Mohammadi and J. Rauh, {\it Prime splittings of determinantal ideals.} Comm. Algebra 46 (2018), no. 5, 2278--2296.

\bibitem{S+Z} B. Sturmfels, A. Zelevinsky, {\it Maximal minors and their leading terms.} Adv. Math. 98 (1993), no. 1, 65--112.

\end{thebibliography}
\end{document}